\documentclass[12pt]{amsart}
\usepackage{amsmath,amssymb,amsthm,amscd}
\numberwithin{equation}{section}

\DeclareMathOperator{\supp}{Supp}

\def\C{{\mathfrak C}}

\def\G{{\mathfrak G}}
\def\g{\mathfrak g}

\def\R{{\mathfrak R}}

\newtheorem{theo}{Theorem}[section]
\newtheorem{lemm}{Lemma}[section]
\newtheorem{coro}{Corollary}[section]

\newtheorem{defi}{Definition}[section]

\def\begeq{\begin{equation}}
\def\endeq{\end{equation}}
\def\p{\partial}
\def\C{\mathbb{C}}

\def\lf{\left}
\def\ri{\right}
\def\e{\epsilon}
\def\ol{\overline}
\def\R{\Bbb R}

\def\ii{\sqrt{-1}}
\def\la{\langle}
\def\ra{\rangle}

\def\bi{{\bar i}}
\def\bj{{\bar j}}
\def\bk{{\bar k}}
\def\p{{\partial}}
\def\bl{{\bar l}}
\def\bla{{\bar \lambda}}
\def\bmu{{\bar\mu}}
\def\lam{{\lambda}}
\def\bnabla{{\overline{\nabla}}}
\def\a{{\alpha}}
\def\ba{{\bar\alpha}}
\def\be{{\beta}}
\def\bbe{{\bar\beta}}

\def\vv<#1>{\langle#1\rangle}
\def\G{\Gamma}
\def\g{\gamma}
\def\la{\langle}
\def\ra{\rangle}
\begin{document}
\title{Product of almost-Hermitian manifolds}

\author{Xu-Qian Fan$^1$}
\address{Department of Mathematics, Jinan University, Guangzhou, 510632, China}
\email{txqfan@jnu.edu.cn}
\thanks{$^1$Research partially supported by the National Natural Science Foundation of China (10901072, 11101106)}
\author{Luen-Fai Tam$^2$}
\address{The Institute of Mathematical Sciences and Department of
 Mathematics, The Chinese University of Hong Kong,
Shatin, Hong Kong, China.} \email{lftam@math.cuhk.edu.hk}
\thanks{$^2$Research partially supported by Hong Kong  RGC General Research Fund
\#CUHK 403011}
\author{Chengjie Yu$^3$ }
\address{Department of Mathematics, Shantou University, Shantou, Guangdong, 515063, China}
\email{cjyu@stu.edu.cn}
\thanks{$^3$Research partially supported by the National Natural Science Foundation of
China (11001161) and GDNSF (9451503101004122).}

\renewcommand{\subjclassname}{%
  \textup{2000} Mathematics Subject Classification}
\date{September 2011}

\begin{abstract}
This is a continuous work about the nonexistence of some complete
metrics on the product of two manifolds studied by Tam-Yu \cite{ty}.  Motivated by the result of Tossati \cite{vt}. We generalize the corresponding
results of Tam-Yu \cite{ty} to the almost-Hermitian case.
\end{abstract}
\maketitle\markboth{Xu-Qian Fan,   Luen-Fai
Tam and Chengjie Yu}{Product of almost-Hermitian manifolds}


\section{Introduction}
 In \cite{Yang}, Yang proved the nonexistence of complete K\"ahler metrics with holomorphic bisectional curvature bounded between two negative constants on the polydisc. Later, Seshadri \cite{Se06} and Seshadri-Zheng \cite{sf} extended Yang's result onto the product of two complex manifolds of positive dimensions. Indeed, they showed that there is no complete Hermitian metrics with holomorphic bisectional curvature bounded between two negative constants and bounded torsion on the product of two complex manifolds of positive dimensions. In \cite{ty}, Tam-Yu relaxed the curvature bounds of the result of Seshadri-Zheng \cite{sf} to a reasonable curvature decay or growth rate in the K\"ahler category. In \cite{vt}, Tosatti generalized the result of Seshadri-Zheng \cite{sf} onto the product two almost complex manifolds. Indeed, Tosatti obtained the following result:
\begin{theo}[Tosatti]
Let $M=X\times Y$ be a product of almost complex manifolds of positive dimensions. Then, there is no complete almost Hermitian metric on $M$ satisfying the following conditions:
\begin{enumerate}
\item The holomorphic bisectional is bounded between two negative constants;
\item The (2,0) part of the curvature tensor is bounded;
\item The torsion is bounded.
\end{enumerate}
\end{theo}

In this paper, motivated by the result of Tosatti, we generalize the results of Tam-Yu \cite{ty} onto the product of two almost complex manifolds of positive dimensions. The main results we obtain are the follows.
\begin{theo}\label{pm1}
Let $X^{2m},Y^{2n}$ be two almost complex manifolds of real dimension
$2m,\ 2n$ respectively, $m,\ n\geq 1$. Then  there is no complete
almost Hermitian metric on
$M=X\times Y$ satisfying the following:
\begin{enumerate}
\item second Ricci curvature$\geq -A(1+r)^2$;
\item holomorphic bisectional curvature$\leq -B<0$;
\item torsion bounded by $A(1+r)$;
\item (2,0) part of the curvature tensor bounded by $A(1+r)^2$
\end{enumerate}
where $r(x)=d(x,o)$ is the distance between $x$ and a fixed point $o\in
M$, and $A, B$ are two positive constants.
\end{theo}

\begin{theo}\label{pm2}
Let $M=X^{2m}\times Y^{2n}$ be the product of two almost complex
manifolds with
positive dimension. Then there is no complete almost Hermitian metric
on $M$ satisfying the following:
\begin{enumerate}
\item second Ricci curvature $\geq -A(1+r^2)^\gamma$;
\item holomorphic bisectional curvature $\leq -B(1+r^2)^{-\delta}$;
\item nonpositive
sectional curvature for the Levi-Civita connection;
\item torsion is bounded by $A(1+r^2)^{\gamma/2}$;
\item (2,0) part of the curvature tensor is bounded by
    $A(1+r^2)^\gamma$
\end{enumerate}
where $\gamma\geq 0 $, $\delta>0$
such that $\gamma+2\delta<1$, $A,B$ are some positive constants, and
$r(x)=d(x,o)$ is the distance of $x$ and a fixed point $o\in M$.
\end{theo}

Clearly, if $(M,J,g)$ is K\"ahler, then the torsion and (2,0) part of
the curvature tensor are zero, and second Ricci curvature is just the
Ricci curvature, so these theorems cover Theorem 1.2 and Theorem 1.3
in \cite{ty} respectively.

The techniques using for proving the main results are mainly the same as in Tam-Yu \cite{ty}. However, for the almost Hermitian case, we don't have at hand a simple formula similar to the K\"ahler or Hermitian case for computing the curvature tensor, so the computation in \cite{ty} can not be extended directly to the almost Hermitian case. In this paper, we use a general Ricci identity (Lemma \ref{lem-Ricci identity}) to handle this difficulty. Generally speaking, this is a Bochner technique on almost Hermitian manifolds. Another difference with the complex case is that we don't have holomorphic coordinates in the almost Hermitian case so that computations can be performed on the coordinate since the complex structure may not be integrable. In this paper, we introduce  local coordinates that play  similar roles of holomorphic coordinates on almost complex manifolds so that similar computations as in K\"ahler geometry can also be performed on almost Hermitian manifolds.

The contents of this paper are arranged as follows. In section 2, we recall some preliminary definitions and results about almost Hermitian manifolds. In section 3, we give a proof of Theorem \ref{pm1}. In section 4, we give a proof of Theorem \ref{pm2}.
\section{Preliminaries  on almost Hermitian manifolds}
For convenience, let us recall some notations and basic results about
almost-Hermitian manifolds,  please see e.g.  \cite{g,twy,vt}.

We say that  $(M^{2n}, J, g)$ is an {\it almost-Hermitian manifold} of real
dimension $2n$ if $J$ is an  almost complex structure on $M$ and $g$ is
a Riemannian metric which is $J$ invariant. For a point $p\in M$,
let  $T_p^\mathbb{C}M=T_pM\otimes\mathbb{C}$, and decompose it as
$T_p^\mathbb{C}M=T'_pM\oplus T''_pM$ where $T'_pM$  and $T''_pM$ are
the eigenspaces of $J$ corresponding to the eigenvalues $\sqrt{-1}$ and
$-\sqrt{-1}$ respectively.

An affine connection $\nabla$ on $TM$ which is extended linearly to
$T^\mathbb{C}M$ is called an {\it almost-Hermitian connection} if
$\nabla J=\nabla g=0.$ Let $\tau$ be the torsion of the connection
$\nabla$ which is defined by
$$\tau(X,Y)=\nabla_XY-\nabla_YX-[X,Y]$$ for $X,Y\in T^\mathbb{C}M$.
One has the following result (see, e.g. \cite{g,k2}).

\begin{lemm}
There exists a unique almost-Hermitian connection $\nabla$ on $(M,J,g)$
such that the torsion $\tau$ has vanishing $(1,1)$ part.
\end{lemm}
This connection is called the {\it canonical connection}. It is first
introduced by Ehresmann and Libermann in \cite{e}, and if $J$ is
integrable it is the connection defined in \cite{c} by Chern. In this
work, we always denote the canonical connection by $\nabla$ and the
Levi-Civita connection by $D$. For the difference between the canonical connection and the Levi-Civita connection, we have the following conclusion, see \cite{g}.
\begin{lemm}\label{lem-comp-connection} On an almost Hermitian manifold
$(M, J, g)$:
$$\vv<D_YX-\nabla_Y
X,Z>=\frac{1}{2}[\vv<\tau(X,Y),Z>+\vv<\tau(Y,Z),X>-\vv<\tau(Z,X),Y>]$$
\end{lemm}
\begin{proof}
Note that
\begin{equation}
\vv<D_XY,Z>+\vv<D_XZ,Y>=X\vv<Y,Z>=\vv<\nabla_XY,Z>+\vv<\nabla_XZ,Y>
\end{equation}
Hence
\begin{equation}\label{eqn-XYZ}
\vv<D_XY,Z>+\vv<D_ZX,Y>=\vv<\nabla_XY,Z>+\vv<\nabla_ZX,Y>-\vv<\tau(Z,X),Y>
\end{equation}
Similarly, we have
\begin{equation}\label{eqn-ZXY}
\vv<D_ZX,Y>+\vv<D_YZ,X>=\vv<\nabla_ZX,Y>+\vv<\nabla_YZ,X>-\vv<\tau(Y,Z),X>
\end{equation}
and
\begin{equation}\label{eqn-YZX}
\vv<D_YZ,X>+\vv<D_XY,Z>=\vv<\nabla_YZ,X>+\vv<\nabla_XY,Z>-\vv<\tau(X,Y),Z>
\end{equation}
Adding (\ref{eqn-XYZ}) and (\ref{eqn-YZX}), and subtracting
(\ref{eqn-ZXY}), we get
\begin{equation}
\vv<D_XY-\nabla_XY,Z>=\frac{1}{2}(\vv<\tau(Y,X),Z>+\vv<\tau(X,Z),Y>-\vv<\tau(Z,Y),X>).
\end{equation}
This completes the proof.
\end{proof}
In   local frame  $e_a$, $1\le a\le 2n$, we have
\begin{coro}
$(\g_{ab}^e-\G_{ab}^e)g_{ec}=\frac{1}{2}(\tau_{ab}^{e}g_{ec}+\tau_{bc}^eg_{ea}-\tau_{ca}^eg_{eb})$
where $\g_{ab}^c$'s are the Christofel symbol of the Riemannian
connection and $$\nabla_{e_a}e_b=\G_{ab}^c e_c.
$$
\end{coro}

\begin{coro}\label{cor-hessian} Let $f$ be a smooth function on $M$,
then
\begin{equation}
\nabla^2f(X,Y)-D^2f(X,Y)=\frac{1}{2}[\vv<\tau(X,Y),\nabla
f>+\vv<\tau(Y,\nabla f),X>-\vv<\tau(\nabla f,X),Y>]
\end{equation}
In local frame $e_a$,

\begin{equation}
f_{ab}-f_{;ab}=\frac{1}{2}(\tau_{ab}^cf_c+\tau_{bc}^e{g^{cd}}f_d
g_{ea}-\tau_{c
a}^eg^{cd}f_dg_{eb})
\end{equation}
where ``$;$" means taking covariant derivatives with respect to the
Levi-Civita connection.

If $\{e_1,\cdots, e_n\}$ is a local frames of $T'M$, then
\begin{equation}
f_{i\bj}-f_{;i\bj}=\frac{1}{2}(\tau_{i\lam}^kg^{\bmu\lam}f_\bmu
g_{k\bj}+\tau_{\bj\bla}^\bk g^{\bla\mu}f_\mu g_{i\bk}).
\end{equation}
In particular, $f_{ab}-f_{ba}=\tau_{ab}^c f_c $ and $f_{i\bj}=f_{\bj
i}$.
\end{coro}
The last assertion follows from the fact that the (1,1) part of $\tau$
is zero.

Taking trace of the above gives us, see \cite{vt}.
\begin{coro}\label{cor-Laplacian}
$\Delta f-\Delta^L f=\tau_{ab}^ag^{bc}f_c$, where $\Delta
f=g^{ab}\nabla^2f(e_a,e_b)$ and $\Delta^Lf=g^{ab}D^2f(e_a,e_b)$ is the
Laplacian with respect to the Levi-Civita connection.
\end{coro}

In general, let $M$ be a manifold with connection $\nabla$, and $E$ be
a vector bundle over $M$ with connection $D$. Let $s$ be a section of
$E$. Then $Ds$ is a section of $TM\otimes E$. To compute more
derivatives, we need the connection $\nabla$ on $M$. Let $\tau$ be the
torsion of $\nabla$. Then, we have following Ricci identity.
\begin{lemm}\label{lem-Ricci identity}
$D^2s(X,Y)-D^2s(Y,X)=-R(X,Y)s+D_{\tau(X,Y)}s$
\end{lemm}
\begin{proof}
\begin{equation}
\begin{split}
&D^2s(X,Y)-D^2s(Y,X)\\
=&(D_YDs)(X)-(D_XDs)(Y)\\
=&D_YD_Xs-Ds(\nabla_YX)-D_XD_Ys+Ds(\nabla_XY)\\
=&D_YD_Xs-D_XD_Ys-D_{[Y,X]}s+Ds(\tau(X,Y))\\
=&-R(X,Y)s+D_{\tau(X,Y)}s
\end{split}
\end{equation}
\end{proof}

Now let us recall some definitions about the curvature. At a point $p$,
choose a local unitary frame $\{e_1,\cdots,e_n\}$ for $T'_p(M)$, and
denote $\{\theta^1,\cdots,\theta^n\}$ as a dual coframe.
Denote
$$R(X,Y)Z=\nabla_X\nabla_Y Z-\nabla_Y\nabla_XZ-\nabla_{[X,Y]}Z,\quad
R(e_C,e_D)e_A=R_{A\ CD}^{\ E}e_E$$
and $R_{ABCD}=R(e_A,e_B,e_C,e_D)=\vv<R(e_C,e_D)e_A,e_B>=R_{A\ CD}^{\
E}g_{EB}.$
Here $A, B, C, D$ can be taken $1,\bar 1,\cdots, n,\bar n.$
Define the {\it second Ricci curvature} as $R'_{k\bar l}=R^{\ l}_{k\
i\bar i}$, the {\it holomorphic bisectional curvature} in the
directions $X$ and $Y$ as
$$B(X,Y)=\frac{R(X,\overline{X},Y,\overline{Y})}{\|X\|^2\|Y\|^2},$$
and the $(2,0)$ part of the curvature as $R^{\ i}_{j\
kl}\theta^k\wedge\theta^l.$

Similar to \cite{vt}, we say that the holomorphic bisectional curvature
is bounded from above by $K$ if
$$B(X,Y)\leq K$$
for all $X,\ Y\in T'M$. The second Ricci curvature is bounded from below by
$-A_1$ if
$$R'_{k\bar l}X^k\overline{X^l}\geq -A_1\|X\|^2$$
for all $X\in T'M.$ The torsion is bounded by $A_2>0$ if
$$|\tau(X,Y)|\leq A_2\|X\|\|Y\|$$
for all $X,\ Y\in T'M.$ The $(2,0)$ part of the curvature is bounded by
$A_3>0$ if
$$|R(\overline{X},Y,Y,X)|\leq A_3\|X\|^2\|Y\|^2$$
for all $X,\ Y\in T'M.$

\section{Proof of Theorem \ref{pm1}}
We will prove Theorem \ref{pm1} by contradiction as in \cite{ty}.
Suppose $M=X^{2n}\times Y^{2m}$ is a product of two almost complex manifolds of positive dimensions satisfying the
conditions in  Theorem \ref{pm1}. Fix a point $q\in Y$, we will show that the volume growth of $X\times\{q\}$ has some upper estimate. On the
other hand we show that this upper estimate is not possible because of
the following maximum principle which is similar to Theorem 1.1
in \cite{kt}.
\begin{lemm} \label{thm-max principle}
Let $(M,g)$ be a complete non-compact Riemannian manifold,  $r(x)$ be
the distance function from a fixed point $o\in M$. Let $u$ be a smooth
function on $M$ satisfying the inequality
\begin{equation}\label{eq:lap-u}
\Delta^L u\geq C_1u^2-C_2(1+r)|\nabla u|
\end{equation}
on $\{u>\delta\}\neq \emptyset$ for some $C_1,  C_2,  \delta>0,$ where
$\Delta^L$ is the Laplace operator with respect to the Levi-Civita
connection, then
\begin{equation} \label{t1-l1}
\liminf_{t\to+\infty}\frac{\log V_o(t)}{t^2}=+\infty
\end{equation}
where $V_o(t)$ is the volume of the geodesic ball of radius $t$ centered
at the point $o\in M$.
\end{lemm}

\begin{proof}
We will adapt the proof of Theorem 2.1 in \cite{kt}. For simplicity, in
the proof of this lemma, we write $\Delta$ instead of $\Delta^L$.
Firstly, we may assume that $\sup_Mu=+\infty$ satisfying
 the  differential inequality \eqref{eq:lap-u} with a different $C_1$.
Otherwise, suppose that
$\sup_Mu=u^*$. By differential inequality \eqref{eq:lap-u} satisfied by $u$,   $u^*$ can
not be attained. Let $v=\frac{1}{u^*-u}$,
we have
\begin{equation*}
\begin{split}
\Delta v=\frac{\Delta u}{(u_*-u)^2}+\frac{2|\nabla u|^2}{(u^*-u)^3}\geq
C_1\delta^2v^2-C_2(1+\rho)|\nabla v|
\end{split}
\end{equation*}
on $\{v>1/(u^*-\delta)\}=\{u>\delta\}$.
 Now $\sup_Mu=\infty$, for any number $Q>\delta$, we can assume that
 $\{u>Q\}$ is not empty. Replace $u$ by $u/Q$,
we know that
\begin{equation*}
\Delta u\geq C_1Qu^{2}-C_2(1+\rho)|\nabla u|
\end{equation*}
on $\{u>1\}.$
So  we conclude that, for any constant $\beta>C_1\delta$, there is a
smooth function $u$ on $M$
such that
\begin{equation}\label{pt1-l1}
\Delta u\geq \beta u^2-C_2(1+r)|\nabla u|
\end{equation}
on the nonempty set $M^*=\{u>1\}$. We will choose $\beta$ to be large
enough later. Note that $C_2$ is independent of $\beta$.

 The same as \cite{kt}, let $0\le\sigma\le 1$ be a smooth function on
 $\R$ such that $\sigma=0$ on $t\le 1$, $\sigma=1$ on $t\ge 2$ and
 $\sigma>0$ for $t>1$, $\sigma'\ge 0$. Let
 $$\lambda(t)=\int_{-\infty}^t \sigma(s)ds.$$
 Then
\begin{displaymath}
 \left\{ \begin{array}{ll}
\lambda(t)\equiv 0  & \textrm{if $t\leq 1$}\\
\lambda(t)>0, \lambda'(t)>0, \lambda''(t)\geq 0 & \textrm{if $t> 1$}\\
\lambda'(t)\equiv 1 & \textrm{if $t\ge2$}.
\end{array} \right.
\end{displaymath}
For $\rho>0$, let $\omega$ be  a Lipschitz continuous   on $M$ such
that
\begin{displaymath}
 \left\{ \begin{array}{ll}
&0\leq \omega\leq 1,\ |\nabla \omega|\leq 1/\rho\\
&\supp(\omega)\subset\overline{B_o(2\rho)}\\
&\omega\equiv 1  \textrm{ on } B_o(\rho).
\end{array} \right.
\end{displaymath}
For any positive constants $p, q, \e$, by \eqref{pt1-l1}, we have
\begin{equation}\label{eq-lambda}
\begin{split}
&\textrm{div}(\omega^{2q}\nabla \lambda(u^p)) \\
=&\la\nabla \omega^{2q},\nabla\lambda(u^p)\ra+\omega^{2q}\Delta
\lambda(u^p)\\
=&2qp\lambda'\omega^{2q-1}u^{p-1}\vv<\nabla\omega,\nabla u>\\
&\quad +\omega^{2q}[\lambda''(pu^{p-1})^2|\nabla
u|^2+p\lambda'(p-1)u^{p-2}|\nabla u|^2+\lambda'pu^{p-1}\Delta u]\\
\geq& p\lambda'\bigg[-\e\omega^{2q}u^{p-2}|\nabla
u|^2-\frac{q^2}{\e}\omega^{2(q-1)}u^p|\nabla
\omega|^2+\omega^{2q}(p-1)u^{p-2}|\nabla u|^2\\
&+\beta
\omega^{2q}u^{p+1}-\omega^{2q}u^p-C_2^2(1+\rho)^2\omega^{2q}u^{p-2}|\nabla
u|^2\bigg]\\
=&p\lambda'\bigg[(\beta-1)\omega^{2q}u^{p+1}+\lf(p-1-C_2^2(1+\rho)^2-\e\ri)\omega^{2q}u^{p-2}|\nabla
u|^2-\frac{q^2}{\e}\omega^{2(q-1)}u^p|\nabla \omega|^2\bigg]
\end{split}
\end{equation}
in $B_o(2\rho)$, provided $u>1$. Since $\lambda'=0$ if $t\le1$ and
$\omega$ has support in $B_o(2\rho)$, the above inequality is true in
$M$.

Choosing $p=p(\rho)$ such that
\begin{equation}\label{pt1-l2}
p-1=2C_2^2(1+\rho)^2,\ \e=\frac{p-1}{2}, q=p+1
\end{equation}
Since $\omega$ has compact support, assume further that $\beta>1$

\begin{equation*}
\begin{split}
\int_{B_o(\rho)}\lambda'\le &\int_{B_o(\rho)} \lambda' (\beta-1)
u^{p+1}\\
\le& \int_{B_o(2\rho)}\lambda' (\beta-1)\omega^{2q}u^{p+1}\\
\end{split}
\end{equation*}
On the other hand, let $q=p+1$, by \eqref{eq-lambda}
\begin{equation*}
\begin{split}
(\beta-1)\int_{B_o(2\rho)}\lambda'
\omega^{2q}u^{p+1}=&(\beta-1)\int_{B_o(2\rho)}\lambda'
\omega^{2(q+1)}u^{p+1}\\
\le &\frac{2(p+1)^2}{p-1}\int_{B_o(2\rho)}\lambda'
\omega^{2p}u^p|\nabla \omega|^2\\
\le& \frac{2(p+1)^2}{ (p-1)\rho^2}\lf(\int_{B_o(2\rho)}\lambda'
\omega^{2(p+1)}u^{p+1}\ri)^{\frac{p}{p+1}}\lf(\int_{B_o(2\rho)}\lambda'\ri)^{\frac1{p+1}}.
\end{split}
\end{equation*}
Hence
\begin{equation*}
    \int_{B_o(\rho)}\lambda'\le
    \lf(\frac{2(p+1)^2}{(\beta-1)(p-1)\rho^2}\ri)^p\cdot\frac{2(p+1)^2}{(p-1)\rho^2}\int_{B_o(2\rho)}\lambda'.
\end{equation*}
By the definition of $p=p(\rho)$, choose $\beta$ such that
$\beta-1>16(C_2^2+1)^2$. There is $\rho_0$ such that if $\rho\ge
\rho_0$,
$$
\frac{2(p+1)^2}{(\beta-1)(p-1)\rho^2}<\frac1{2(C_2^2+1)}.
$$
Hence for $\rho\ge \rho_0$,
$$
\int_{B_o(\rho)}\lambda'(u^{p(\rho)}) \le \lf(\frac12\ri)^{p}
\int_{B_o(2\rho)}\lambda'(u^{p(\rho)})\le
\int_{B_o(2\rho)}\lambda'(u^{p(2\rho)}).
$$
for some $k$ which is independent of $\rho$. Here we have used the fact
that $\lambda'$ is nondecreasing, $p(2\rho)>p(\rho)$ and $\lambda'=0$
if $u\le 1$. Let
$$
F(\rho)=\int_{B_o(\rho)}\lambda'(u^{p(\rho)}).
$$
We have
$$
F(\rho)\le \lf(\frac12\ri)^{p} F(2\rho).
$$
By iterating, we have
$$
F(\rho_0)\le \lf(\frac12\ri)^{C_3\rho^2}F(\rho)\le
\lf(\frac12\ri)^{C_3\rho^2}V_o(\rho)
$$
for some $C_3>0$, for $\rho>0$, because $\lambda'\le 1$. Since
$\{u>1\}$ is nonempty, if $\rho_0$ is chosen large enough,
$F(\rho_0)>0$. From this it is easy to see that the lemma is true.
\end{proof}

To estimate the volume growth, we will use the following result   due
to Tosatti \cite[Theorem 4.2]{vt}:
\begin{lemm}[Tossati]\label{lcle1}
Suppose $(M^m,J,g)$ is a complete almost-Hermitian manifolds with real
dimension $m$, $B_o(R)$ is a geodesic ball centered at $o\in M$ of
radius $R$. If the second Ricci curvature of $B_o(R)$ is bounded from
below by $-K_1$,  the torsion bounded by $A_2$ and $(2,0)$ part of the
curvature bounded by $A_3$ on $B_o(R)$ for some positive constants
$A_1, A_2, K$, then
\begin{equation} \label{lemlc1}
\Delta r\leq \frac{m}{r}+c
\end{equation}
where  $r$ is   the distance function from $o$,  $\Delta$ is the
Laplace operator with respect to the canonical connection,
$c=c_1\alpha$, $\alpha=(A_2+\sqrt{K_1}+\sqrt{A_3})$, and $c_1$ is a
positive constant  depending on $m$.
\end{lemm}

From this, one has
\begin{coro} \label{corvc1}
Under the same notations and assumptions as in Lemma \ref{lcle1}, we
have for any fixed $0<t_0<R$,
\begin{equation*}
  V_o(t)\leq V_o(t_0)\left(\frac{t}{t_0}\right)^{m+1}e^{C\alpha t}
  \textrm{ for } R\geq t\geq t_0,
\end{equation*}
where $C$ is a constant depending only on $m$. In particular,
\begin{equation*}
  V_o(R)\leq V_o(t_0)\left(\frac{R}{t_0}\right)^{m+1}e^{C\alpha R}.
\end{equation*}
\end{coro}

\begin{proof}
By Corollary \ref{cor-Laplacian} or Lemma 3.2 in \cite{vt}, and by
Lemma \ref{lcle1}, we have
$$\Delta r\leq \frac{m}{r}+c $$
where $c$ is the constant in Lemma \ref{lcle1}.
So
\begin{equation*}
  \Delta^L r\leq \frac{m}{r}+C\alpha
\end{equation*}
where $C$ is a positive constant depending only on $m$. Multiplying $r$ to the both
sides of inequality above, we have
\begin{equation*}
  r\Delta^L r\leq m+r C\alpha.
\end{equation*}
So
\begin{equation*}
  \int_{B_o(t)}r\Delta^L r\leq \int_{B_o(t)}\left(m+r C\alpha\right).
\end{equation*}
Hence
$$tA_o(t)\leq \left(m+1+ C\alpha t \right)V_o(t).$$
That is
$$(\ln V_o(t))'\leq t^{-1}\left(m+1+C\alpha t \right).$$
Integrating both sides from $t_0$ to $r$, we have
$$V_o(t)\leq V_o(t_0)\left(\frac{t}{t_0}\right)^{m+1}e^{C\alpha t}.$$
This completes the proof of the lemma.
\end{proof}

Similar to Lemma 2.1 in \cite{ty}, we have
\begin{lemm} \label{lemls1}
Suppose $(M^m,J,g)$ and $(N^n,\tilde{J},h)$ are two complete
almost-Hermitian manifolds. Let $f$ be a non-constant almost-complex
map from $M$ to $N$. Let $o\in M$ and let $R>0$.  If the second Ricci
curvature of $B_o(2R)$ is bounded from below by $-K_1$,  the torsion
bounded by $A_2$ and $(2,0)$ part of the curvature bounded by $A_3$ on
$B_o(2R)$, and the  bisectional curvature in $f(B_o(2R))$ is bounded
above by $-K_2$, where $A_2, A_3, K_1,\ K_2$ are positive constants,
then on $B_o(R)$,
\begin{equation*}
  f^*h\leq \frac{2K_1+CR^{-2}(1+cR)}{2K_2}g
\end{equation*}
where $C$ is a constant depending only on $m$ and $c$ is the same as in
\eqref{lemlc1}.
\end{lemm}

\noindent
\begin{proof}
Let $u=\textrm{tr}_g(f^*h)$. Since  the second Ricci curvature of
$B_o(2R)$ is bounded from below by $-K_1$ and the bisectional curvature
in $f(B_o(2R))$ is bounded above by $-K_2$, by the result of \cite[page
1081]{vt}, one has
\begin{equation}\label{lsch1}
  \triangle u\geq 2K_2u^2-2K_1u
\end{equation}
on $B_o(2R)$, where $\Delta$ denotes the Laplacian with respect to the canonical
connection on $M$.

Similar to the proof of Lemma 2.1 in \cite{ty}. Let $\eta\geq 0$ be a
smooth function on $\mathbb{R}$ such that (1) $\eta(t)=1$ for $t\leq
1$, (2) $\eta(t)=0$ for $t\geq 2$, (3) $-C_3\leq \eta'/\eta^{1/2}\leq
0$ for all $t\in \mathbb{R}$, and (4) $|\eta''(t)|\leq C_3$ for all
$t\in \mathbb{R}$ for some absolute constant $C_3>0$. Let
$\phi=\eta(r/R)$, where $r$ is the distance function from $o$.

Suppose $\phi u$ attains maximum at $\bar{x}\in B_o(2R)$, then
$\phi(\bar{x})>0.$ Using an argument of Calabi as in \cite{cy}, we may
assume that $\phi u$ is smooth at $\bar{x}$. Then we have (1)
$\nabla(\phi u)(\bar{x})=0$ which implies that at $\bar{x}$, $\nabla
u=-u\phi^{-1}\nabla \phi,$ (2) $\Delta (\phi u)(\bar{x})\leq 0$. Using
Corollary \ref{cor-Laplacian},   at $\bar{x}$, we have
\begin{align}
0&\geq \Delta (\phi u) \nonumber\\
&= \phi \Delta u+u\Delta \phi+2\vv<\nabla\phi,\nabla u>  \\
&=\phi \Delta u+u(\eta''R^{-2}+\eta'R^{-1}\Delta
r)+2\vv<\nabla\phi,\nabla u> \nonumber\\
&=\phi \Delta u-2uR^{-2} \frac{(\eta')^2}{\eta}
+u(\eta''R^{-2}+\eta'R^{-1}\Delta r)\nonumber\\
&\geq \phi \Delta u-2C_3^2uR^{-2}+u(\eta''R^{-2}+\eta'R^{-1}\Delta r
)\nonumber\\
&\geq \phi (2K_2u^2-2K_1u)-2C_3^2uR^{-2}-C_3uR^{-2}+u\eta'R^{-1}\Delta
r  (\textrm{by} \eqref{lsch1}).\nonumber
\end{align}
So
$$ 2\phi K_2 u^2 \leq 2K_1\phi
u+2C_3^2uR^{-2}+C_3uR^{-2}-u\eta'R^{-1}\Delta r.$$
By Lemma \ref{lcle1}, we have
$$\Delta r\leq \frac{m}{r}+c$$
where $c$ is the same as in \eqref{lemlc1}.
So  we can get
\begin{align}
 2K_2\phi u^2&\leq 2K_1\phi
 u+2C_3^2uR^{-2}+C_3uR^{-2}+C_3R^{-1}u\left(\frac{m}{R}+c\right)
 \nonumber\\
 &\leq u[2K_1+R^{-2}(2C_3^2+C_3+C_3m+C_3cR)].
 \end{align}
 Hence
 \begin{align}
 \sup_{B_o(R)}u\leq\sup_{B_o(2R)}(\phi u)\leq
 \frac{2K_1+CR^{-2}(1+cR)}{2K_2}
 \end{align}
 where $C$ is a constant depending on $m$ and $C_1$, and $c$ is the
 same as in \eqref{lemlc1}. Therefore the lemma holds.
\end{proof}

 We have the following volume growth estimate of geodesic ball on the
 submanifold $X$.
\begin{lemm}\label{volume-l1}
Let $M, X, Y$ as in Theorem \ref{pm1}. Suppose  there is a complete
almost Hermitian metric on
$X\times Y$ satisfying the assumptions in the theorem, that is:
\begin{enumerate}
\item second Ricci curvature$\geq -A(1+r)^2$
\item holomorphic bisectional curvature$\leq -B<0$
\item torsion bounded by $A(1+r)$
\item (2,0) part of the curvature tensor bounded by $A(1+r)^2$
\end{enumerate}
for some positive constants $A, B$,
where $r(x)=d(x,o)$ is the distance of $x$ and a fixed point
$o=(p,q)\in X\times Y$ . Let $V^{^{X_q}}_p(t)$ be the volume of the
geodesic ball of radius $t$ with center at $p$ with respect to the
induced metric $g^q$ on $X_q=X\times \{q\}$. Then
$$
V^{^{X_q}}_p(t)\le C_4\exp(C_4 t^2)
$$
for some positive constant $C_4$  independent of $t$.
 \end{lemm}
\begin{proof}
The proof is similar to the proof of Lemma 2.2 in \cite{ty}. By \cite{k},
the bisectional curvature of the canonical connection
 of an almost complex submanifold is not bigger than the one of the
 ambient space.  So for any point $y_0\in Y$, then the bisectional
 curvature of $X_{y_0}$ is also bounded from above by $-B$.  The same
 holds for $Y_{x_0}$ where $Y_{x_0}=\{x_0\}\times Y$. By Lemma
 \ref{lemls1},
 there is a constant $C_5$ independent of $x_0, y_0$ such that
\begin{equation*}
\begin{split}
(\pi''_{x_0})^*(g^{ x_0})|_{(x,y)}&\le C_5\lf(1+r
(x,y)\ri)^2g|_{(x,y)}\quad \textrm{ and}\\
(\pi'_{y_0})^*(g^{y_0})|_{(x,y)}&\le C_5\lf(1+r (x,y)\ri)^2g|_{(x,y)}
\end{split}
\end{equation*}
for  $(x,y)\in M$, where  $\pi'_{y_0}$ is the projection from $M$ onto
$X_{y_0}$ defined by $\pi'_{y_0}(x,y)=(x,y_0)$, $\pi''_{x_0}$ is the
projection from $M$ onto $Y_{x_0}$ defined by
$\pi''_{x_0}(x,y)=(x_0,y)$, $g^{x_0}, g^{y_0}$ are the induced metrics
on $Y_{x_0},\ X_{y_0}$ respectively. By Corollary \ref{corvc1}, we have
$V_o(2R)\leq \exp(C(1+R)^2)$ for some constant $C$. In the rest of the
proof, we can follow the corresponding argument of the proof of Lemma
2.2 in \cite{ty} to get the conclusion of this lemma.
\end{proof}

Next we want to find a function on $X_q$ satisfying the differential
inequality in   Lemma \ref{thm-max principle}. For convenience, we will
introduce a special coordinate near a point on almost complex
manifold.

\begin{defi}
Let $(M^{2m},J)$ be an almost complex manifold. Let $p\in M$, $U$ be an
open neighborhood of $p$. Let $\phi:U\to \Omega\subset \mathbb{C}^m$ be
a diffeomorphism. Then, $(U,\phi)$ is called a complex coordinate.
\end{defi}
Let $\phi=(z^1,z^2,\cdots,z^m)$ be a complex coordinate, and suppose
that $z^i=x^i+\ii y^i$. Then $(x^1,y^1,\cdots,x^m,y^m)$ is a local
coordinate. As usual, we define
\begin{equation}
\frac{\p}{\p z^i}=\frac{1}{2}(\frac{\p}{\p x^i}-\ii\frac{\p}{\p y^i})
\end{equation}
and
\begin{equation}
\frac{\p}{\p z^\bi}=\frac{1}{2}(\frac{\p}{\p x^i}+\ii\frac{\p}{\p
y^i})
\end{equation}
as vectors in $TM\otimes\C$.
\begin{defi}
Let $(M,J)$ be an almost complex manifold. Let $(z^1,z^2,\cdots,z^m)$
be a complex coordinate at $p$. It is called {\it almost holomorphic}
at $p$ if
\begin{equation}
J(\frac{\p}{\p z^i})(p)=\ii\frac{\p}{\p z^i}(p)
\end{equation}
for all $i=1,2,\cdots,m$.
\end{defi}
\begin{lemm}
Let $(M,J)$ be an almost complex manifold. Then for any $p\in M$, there
is a local complex
coordinate $(z^1,z^2,\cdots ,z^n)$ that is almost holomorphic at $p$
such that
\begin{equation}
\p_i J_j^k(p)=\p_\bi J_j^k(p)=\p_\bi J_j^\bk(p)=0
\end{equation}
where
\begin{equation}
J(\frac{\p}{\p z^i})=J_i^j\frac{\p}{\p z^j}+J_i^\bj\frac{\p}{\p
z^\bj}.
\end{equation}
\end{lemm}
\begin{proof}
Let $(z^1,z^2,\cdots,z^n)$ be a local complex
coordinate that is almost holomorphic at $p$. Suppose that
\begin{equation}
J(\frac{\p}{\p z^i})=J_i^{j}\frac{\p}{\p z^j}+J_i^{\bj}\frac{\p}{\p
z^\bj}
\end{equation}
Then
\begin{equation}
J_i^{j}(p)=\ii\delta_{ij}\ \mbox{and}\ J_i^{\bj}(p)=0
\end{equation}
By that $J^2=-id$, we know that
\begin{equation}\label{eqn-J-1}
J_i^{j}J_{j}^k+J_i^{\bj}J_{\bj}^k=-\delta_{ik}.
\end{equation}
Taking partial differentiations of \eqref{eqn-J-1}, we know that
\begin{equation}\label{eqn-first-almost-holo}
\p_j J_i^{k}(p)=\p_{\bj}J_{i}^k(p)=0
\end{equation}
for all $i,j,k=1,2,\cdots,m$.

Let $(w^1,w^2,\cdots,w^n)$ be a coordinate change of
$(z^1,z^2,\cdots,z^n)$ such that
\begin{equation}\label{eqn-first-change}
\frac{\p w^i}{\p z^j}(0)=\delta_{ij}\ \mbox{and}\ \frac{\p w^\bi}{\p
z^j}(0)=0.
\end{equation}
Suppose that
\begin{equation}\label{eqn-first-derivatives}
J(\frac{\p}{\p w^i})=\hat J_i^{j}\frac{\p}{\p w^j}+\hat
J_i^{\bj}\frac{\p}{\p w^\bj}.
\end{equation}
By a straight forward computation, we have
\begin{equation}
\hat J_i^{\bl}=\frac{\p z^j}{\p w^i}J_j^k\frac{\p w^\bl}{\p
z^k}+\frac{\p z^\bj}{\p w^i}J_{\bj}^k\frac{\p w^\bl}{\p z^k}+\frac{\p
z^j}{\p w^i}J_j^\bk\frac{\p w^\bl}{\p z^\bk}+\frac{\p z^\bj}{\p
w^i}J_\bj^\bk\frac{\p w^\bl}{\p z^\bk}.
\end{equation}
Then, by using (\ref{eqn-first-almost-holo}) and
(\ref{eqn-first-change}), we have
\begin{equation}
\begin{split}
&\frac{\p}{\p w^\ba}\hat J_i^\bl(p)\\
=&\ii \frac{\p}{\p w^\ba}(\frac{\p w^\bl}{\p z^i})-\ii\frac{\p}{\p
w^\ba}(\frac{\p z^\bl}{\p w^i})+\frac{\p}{\p w^\ba}(J_i^\bl)\\
=&-2\ii\frac{\p z^\bl}{\p w^\ba\p w^i}(0)+\frac{\p}{\p
z^\ba}(J_i^\bl)(p)
\end{split}
\end{equation}
So, if we choose $(w^1,w^2,\cdots,w^n)$ such that
\begin{equation}\label{eqn-second-derivates}
\frac{\p z^\bl}{\p w^\ba\p w^i}(0)=\frac{1}{2\ii}\frac{\p}{\p
z^\ba}(J_i^\bl)(p)
\end{equation}
and (\ref{eqn-first-derivatives}) are both true. Then
\begin{equation}
\frac{\p}{\p w^\ba}\hat J_i^\bl(p)=0
\end{equation}
and $(w^1,w^2,\cdots,w^m)$ is a complex coordinate that is almost holomorphic
at $p$. This completes the proof.
\end{proof}
\begin{defi}
We call the local coordinate in the last lemma an  holomorphic
coordinate at $p$.
\end{defi}
\begin{coro}\label{cor:prod-coor}
Let $M=X\times Y$ be a product of two almost complex manifolds. Let $(z^1,z^2,\cdots,z^k)$ be a local holomorphic coordinate for $X$ at $x$ and $(w^{1},w^{2},\cdots,w^l)$ be a local holomorphic coordinate for $Y$ at $y$. Then $(z^1,z^2,\cdots,z^k,w^1,\cdots,w^l)$ is a local holomorphic coordinate at $p=(x,y)$.
\end{coro}
\begin{proof}
Let $J_X$ and $J_Y$ be the almost complex structures on $X$ and $Y$ respectively. Since the almost complex structure on $M=X\times Y$ is
a product of $J_X$ and $J_Y$, we have
\begin{equation}
J(\frac{\p}{\p z^i})=J_X(\frac{\p}{\p z^i})=(J_X)_i^j(z)\frac{\p}{\p z^j}+(J_X)_i^\bj(z)\frac{\p}{\p z^\bj}
\end{equation}
and
\begin{equation}
J(\frac{\p}{\p w^\a})=J_Y(\frac{\p}{\p w^\a})=(J_Y)_\a^\be(w)\frac{\p}{\p w^\be}+(J_Y)_\a^\bbe(w)\frac{\p}{\p w^\bbe}.
\end{equation}
Then the conclusion comes directly by a simple computation.
\end{proof}

For almost-Hermitian manifold with canonical connection, we have

\noindent
\begin{lemm}\label{lem-holo coord}
Let $(M^{2n},J,g)$ be an almost Hermitian  complex manifold and
$\nabla$ be the canonical connection. Then, for each point $p\in M$ and any   holomorphic
      coordinate $(z^1,z^2,\cdots,z^n)$ at $p$, we have
\begin{equation}
\nabla_{\frac{\p}{\p \bar z^i}}\frac{\p}{\p z^j}(p)=0.
\end{equation}

\end{lemm}
\begin{proof}  Since $\nabla J=0$,
\begin{equation} \label{spt1-1}
\begin{split}
&J\left(\nabla_{\frac{\p}{\p\bar  z^i}}\frac{\p}{\p z^j}(p)\right)\\
=&\nabla_{\frac{\p}{\p\bar z^i}}\left(J\frac{\p}{\p z^j}\right)(p)\\
=&\nabla_{\frac{\p}{\p\bar z^i}} \left(J_j^k\frac{\p}{\p
z^k}+J_j^\bk\frac{\p}{\p z^\bk} \right)(p)\\
=&J_j^k\nabla_{\frac{\p}{\p \bar z^i}}\frac{\p}{\p z^k}(p)+J_j^\bk
\nabla_{\frac{\p}{\p \bar z^i}}\frac{\p}{\p z^\bk}(p),\text{\ (by the
definition of holomorphic coordinates)}\\
=&\ii \nabla_{\frac{\p}{\p \bar z^i}}\frac{\p}{\p z^j}(p)
\end{split}
\end{equation}
because $J_i^k=\ii \delta_{ij}, \ J_i^{\bar j}=0$ at $p$. Hence
 $\nabla_{\frac{\p}{\p z^\bi}}\frac{\p}{\p z^j}(p)$ is a (1,0) vector.
 Similarly, we can show that
$\nabla_{\frac{\p}{\p z^j}}\frac{\p}{\p\bar  z^i}(p)$ is a (0,1)
vector. On the other hand,
\begin{equation}
\nabla_{\frac{\p}{\p z^j}}\frac{\p}{\p\bar z^i}(p)-\nabla_{\frac{\p}{\p
\bar z^i}}\frac{\p}{\p z^j}(p)=\tau\left(\frac{\p}{\p
z^j}(p),\frac{\p}{\p \bar z^i}(p)\right)=0.
\end{equation}
Hence the conclusion follows.

\end{proof}

From the lemma, one can get the following.
\begin{coro}\label{cor-hessian 2}
Let $(z^1,z^2,\cdots,z^n)$ be a holomorphic coordinate at $p$
on an almost Hermitian
manifold $(M^{2n},J,g)$, then
\begin{equation}
u_{i\bj}(p)=\p_i\p_\bj u(p)
\end{equation}
where $u_{i\bj}=\nabla^2u(\p_i,\p_\bj)$ is the complex Hessian with respect to the canonical
connection.
\end{coro}

We also need the following facts on submanifolds. Let $(M,J,g)$ be an
almost Hermitian manifold and $\overline{\nabla}$ the
canonical connection, and $\bar\tau$ be the torsion of
$\overline\nabla$.
Let $N$ be a submanifold of $M$. Define the connection on $N$
\begin{equation}
\nabla_XY=(\bnabla_XY)^\top
\end{equation}
We will also need the following result about the torsion of
submanifold.
\begin{lemm}\label{lem-submfd}
\begin{itemize}
  \item [(a)] $\nabla$ is the canonical connection of the induced
      almost Hermitian manifold
$(N,J,g)$ with torsion
\begin{equation}
\tau(X,Y)=\bar\tau(X,Y)^\top
\end{equation}
for any $X,Y\in TN$.
  \item [(b)] $h(X,\ol Y)=h(\ol Y,X)=0$ for $X,Y\in T'(N)$, where
      $h(U,V)=-\lf(\ol\nabla_{U}V\ri)^\perp$, $U,V\in T^{\C}(N)$.
  \item [(c)] Let $f$ be a smooth function on $M$, then
  $$
  \nabla^2f(X,\ol Y)=\ol\nabla^2 f(X,\ol Y)
  $$
  for $X, Y\in T'(N)$.
\end{itemize}
\end{lemm}
\begin{proof}
(a) By the definitions of the torsion and the connection $\nabla$, for
$U,V\in T^\mathbb{C}N$, we have
\begin{equation*}
\begin{split}
\tau(U,V)&=\nabla_UV-\nabla_VU-[U,V]\\
&=(\bnabla_UV)^\top-(\bnabla_VU)^\top-[U,V]\\
&=\bar{\tau}(U, V)^\top.
\end{split}
\end{equation*}
 Clearly $\nabla$ is also a canonical connection on $N$.

 (b) Noting that $JW\in T^\mathbb{C}N$ for $W\in T^\mathbb{C}N$, we can
 get
\begin{equation}\label{eqlch1}
h(U,JV)=J(h(U,V))
\end{equation}
for $U,V\in T^\mathbb{C}N$. Since $\tau(X,\bar Y)=0$, we have
\begin{equation*}
\begin{split}
h(\bar Y,X)-h(X,\bar Y)
&=(\bnabla_X\bar Y)^\perp-(\bnabla_{\bar Y}X)^\perp\\
&=(\bnabla_X\bar Y-\bnabla_{\bar Y}X)^\perp\\
&=([X,\bar Y]+\tau(X,\bar Y))^\perp=0.
\end{split}
\end{equation*}
So
\begin{equation}\label{eqlch3}
h(X,\bar Y)=h(\bar Y,X).
\end{equation}
Let $\{e_1,\cdots, e_n\}$ be a unitary frame on $T'N$ where $n$ is the
complex dimension of $N$, by \eqref{eqlch1}, we can get
\begin{equation}\label{eqlch4}
\begin{split}
 \vv<h(e_i,e_{\bar{j}}),e_{\bar{k}}>
&=\vv<Jh(e_i,e_{\bar j}),Je_{\bar k}>\\
&=\vv<h(e_i,Je_{\bar j}),Je_{\bar k}>\\
&=-\vv<h(e_i,e_{\bar{j}}),e_{\bar{k}}>
\end{split}
\end{equation}
and by \eqref{eqlch3},
\begin{equation}\label{eqlch5}
\begin{split}
\vv<h(e_i,e_{\bar{j}}),e_k>
&=\vv<h(e_{\bar j},e_i),e_k>\\
&=\vv<Jh(e_{\bar j},e_i),Je_k>\\
&=\vv<h(e_{\bar j},Je_i),Je_k>\\
&=-\vv<h(e_i,e_{\bar{j}}),e_k>,
\end{split}
\end{equation}
so $\vv<h(e_i,e_{\bar{j}}),e_{\bar{k}}>=0=\vv<h(e_i,e_{\bar{j}}),e_k>$,
here $i,j,k\in\{1,\cdots,n\}.$ Hence

\begin{equation}\label{eqlch2}
h(X,\bar Y)=h(\bar Y,X)=0
\end{equation}
for any $X,Y\in T'N$.

(c)
\begin{equation}\label{eqn-com-sub-hess}
\begin{split}
&\nabla^2f(X,\bar Y)\\
=&\bar YX(f)-\nabla_{\bar Y}X(f)=\bar YX(f)-[\bnabla_{\bar Y}X+h(\bar
Y,X)](f)\\
=&(\bnabla^2f)(X,\bar Y).
\end{split}
\end{equation}
Therefore the lemma is true.
\end{proof}
%

Now we are ready to prove Theorem \ref{pm1}.
\begin{proof}[Proof of Theorem \ref{pm1}.]
We proceed by contradiction. Let $g$ be a complete almost Hermitian
metric on $X^{2m}\times Y^{2n}$ satisfying the assumptions,  $\bnabla$
be the canonical connection.

 Denote the fixed point $o$ as $(p,q)\in X\times Y$. Consider the
 inclusion map: $i:X_q\hookrightarrow X\times Y$ defined by
 $i(x)=(x,q)$, and pull back
 the tangent bundle $T(X\times Y)$ by $i$ on $X_q$. Let $u\in T'_q(Y)$,
 we can get a section $V$ of $i^*T(X\times Y)$ on $X_q$ such that
 $V(x)=u$ for all $x\in X_q$.  For simplicity, let $\{e_1,\cdots,
 e_m,e_{m+1},\cdots, e_{m+n}\}$ be a unitary
 frame on $T'M$ such that $\{e_1,\cdots, e_m\}$ is a  frame on $T'X_q$.
 In the rest of this proof, we will take $\alpha\in\{1,\cdots, m\}$ and
 $i,j\in\{1,\cdots,m+n\}$.  Since $u$ is a $(1,0)$ vector, we can write
 $V=V^ie_i$. By the Ricci identity (Lemma \ref{lem-Ricci identity}),  see also e.g. \cite[page
 1075]{vt} and Lemma 3.1 therein, we have
\begin{equation}\label{t1-pf-v1}
\begin{split}
&\frac{1}{2} \Delta_{X_q}\|V\|^{2}\\
=&(V^i\overline{V^i})_{\a\ba}  \\
=&V^i_{\ \a\ba}\overline{V^i}+V^i\overline{V^i}_{\a\ba}+V^i_{\
\a}\overline{V^i}_{\ba}+V^i_{\ \ba}\overline{V^i}_{\a}\\
=&V^i\overline{V^i}_{\a\ba}+V^i_{\
\ba\a}\overline{V^i}-\overline{V^i}R_{j\bi\a\ba}V^j+V^i_{\
\a}\overline{V^i}_{\ba}+V^i_{\ \ba}\overline{V^i}_{\a}\\
=&-\overline{V^i}R_{j\bi\a\ba}V^j+V^i_{\
\a}\overline{V^i}_{\ba}\\
\geq& mB\|V\|^2
\end{split}
\end{equation}
where we have used Corollary \ref{cor:prod-coor} which implies that $V^i_{\ \ba}=0$ and hence $V^i_{\ \ba\a}=0.$

Moreover, by Corollary \ref{cor-Laplacian} (See also Lemma 3.2 in \cite{vt}), and the assumption of the torsion,
we have
\begin{equation}
\Delta_{X_q}^L\|V\|^{2}\geq  2mB\|V\|^2-C(m, A)(1+\rho)\|\nabla \|V\|^2\|,
\end{equation}
here $\rho$ is the distance function from $p$ on $X_q$. Similar to the
proof (2.8) in \cite{ty} using the Schwartz lemma (Lemma \ref{lemls1}), one can get that $\|V\|^2$ is a positive bounded
function.
By Lemma \ref{thm-max principle} and Lemma \ref{volume-l1}, we have a
contradiction because $|V|>0$.
\end{proof}

\section{Proof of Theorem \ref{pm2}}

We need a lemma similar  to Lemma 3.1 in \cite{ty}. Since $M$ may not
be K\"ahler, we need Corollary \ref{cor-hessian 2} and Lemma \ref{lem-submfd} in our
computations.
\begin{lemm} \label{t2-p-l1}
Let $M=X^{2m}\times Y^{2n}$ be the product of two almost complex
manifolds with
positive dimensions. Assume that $M$ is simply connected. Suppose there is a
complete Hermitian metric $g$ on $M$ satisfying the assumptions in
Theorem \ref{pm2}, that is:
\begin{enumerate}
\item the second Ricci curvature $\geq -A(1+r^2)^\gamma$;
\item the holomorphic bisectional curvature $\leq
    -B(1+r^2)^{-\delta}$;
\item
sectional curvature for the Levi-Civita connection is nonpositive;
\item torsion is bounded by $A(1+r^2)^{\gamma/2}$;
\item (2,0) part of the curvature tensor is bounded by
    $A(1+r^2)^\gamma$;
\end{enumerate}
where $\gamma\geq 0 $, $\delta>0$
such that $\gamma+2\delta<1$, $A,B$ are some positive constants, and
$r(x,y)=d(o,(x,y))$ is the distance of $(x,y)\in X\times Y=M$ from a
fixed point $o\in M$. Then
there is a positive constant $C$ depending only on
$m,n,\gamma,\delta,A$ and $B$, such that
\begin{equation}
g|_{(x_0,y)}(u,\bar u)\leq C(1+r^2(x,y))^\gamma(1+r^2(x_0,y))^\delta
g|_{(x,y)}(u,\bar u)
\end{equation}
for any $x_0,x\in X$, $y\in Y$ and $u\in T'_y(Y)$.
\end{lemm}
\begin{proof}
Let $\pi:X\times Y\to \{x_0\}\times Y=Y_{x_0}$ be the natural
projection. We only need to prove that
\begin{equation}
u(x,y)\leq C(1+r^2(x,y))^\gamma(1+r^2(x_0,y))^\delta
\end{equation}
where $u(x,y)$ is the energy density of $\pi$.

By equation (5.9) in \cite{vt}, the assumptions (1) and (2),  and the
fact that the bisectional curvature of the canonical connection
 of an almost complex submanifold is not bigger than the one of the
 ambient space \cite{k}, we have at $(x,y)\in M$,
\begin{equation}\label{pt2lulb1}
\Delta u\geq -2A(1+r^2(x,y))^\gamma u+2B(1+r^2(x_0,y))^{-\delta}u^2.
\end{equation}

Let $(x,y)\in M$. $T_{(x,y)}'(M)=T_x'(M)\oplus T_y'(M)$,  by Lemma
\ref{lem-holo coord} we can choose   a holomorphic coordinate
$(z^1,z^2,\cdots,z^m)$ of $X$ at $x$ and  a
  holomorphic coordinate $(z^{m+1},\cdots,z^{m+n})$ of $Y$ at $y$. Then, by Corollary \ref{cor:prod-coor},
  $(z^1,\dots,z^{m+n})$ is a holomorphic coordinate at $(x,y)$ in $M$.
Note that
\begin{equation}\label{pt2-ue1}
u(x,y)=g_{\a\bbe}(x_0,y)g^{\bbe\a}(x,y)
\end{equation}
where $\a,\be\in\{m+1,\cdots, m+n\}$. Let $f(y)=r(x_0,y)$,  $y\in
Y_{x_0}$. By abusing notations, we also denote  the function $f\circ
\pi$ on $M$ be $f$. Since $M$ is simply connected with nonpositive
Riemannian curvature, $f^2$  is a smooth function.
Then
\begin{equation}\label{eqn-gradient}
\begin{split}
 |\nabla f |^2(x,y)
\le& e(\pi)|\nabla_{Y_{x_0}} f|^2 (y)\\
\leq&u(x,y).
\end{split}
\end{equation}
 Near $x$, choose a frame $\{e_1,\cdots,e_m\}$ on $T'(X)$ with dual
 co-frame $\omega_1,\dots,\omega_m$, and near $y$ choose a frame
 $\{e_{m+1},\cdots,e_{m+n}\}$ on $T'(Y)$ with the dual co-frame
 $\{\omega^{m+1},\cdots,\omega^{m+n}\}$ satisfying that at $y$, $
 \vv<\omega^\alpha,\overline{\omega^\beta}>|_{(x,y)}=\delta_{\alpha\be}.$
 Then $e_1,\dots,e_{m+n}$ is a frame near $(x,y)$ with coframe
 $\omega^1,\dots,\omega^{m+n}$. Moreover, at $(x,y)$, $e_1,\dots, e_m$
 are linear combinations of $\frac{\p}{\p z^1},\dots,\frac{\p}{\p z^m}$
 and $e_{m+1},\dots,e_{m+n}$ are linear combinations of $\frac{\p}{\p
 z^{m+1}},\dots,\frac{\p}{\p z^{m+n}}$. Without loss of generality we
 may assume that $e_a=\frac{\p}{\p z^a}$ for all $a$ at the point $(x,y)$.

  By Corollary \ref{cor-hessian 2}, the fact that $f^2$ is independence
  of $x$, and Corollary \ref{cor-hessian}, we have
\begin{equation} \label{eqn-laplacian1}
\begin{split}
&\Delta f^2|_{(x,y)}\\
=&2g^{\bar ba}(x,y)\partial_{a}\partial_{\bar b}f^2|_{(x,y)}\\
=&2g^{\bar\beta\alpha}(x,y)\partial_{\alpha}\partial_{\bar
\beta}f^2|_{(x,y)}\\
=&2g^{\bar\beta\alpha}(x,y)(\partial_{\alpha}\partial_{\bar
\beta}r^2)|_(x_0,y)\\
=&2g^{\bar\beta\alpha}(x,y)(r^2)_{;\a\bbe}(x_0,y)\\
& +2g^{\bar\beta\alpha}(x,y)\left[\frac{1}{2}\tau_{\alpha t}^b g^{\bk
t}g_{b\bar \beta}\p_{\bar k}(r^2)+\frac{1}{2}\tau_{\bar \beta\bar
t}^{\bar h} g^{\bar t k}g_{\alpha\bar h}\p_k
(r^2)\right]\bigg|_{(x_0,y)}\\
=&2g^{\bar\beta\alpha}(x,y)(r^2)_{;\a\bbe}(x_0,y)\\
& +g^{\bar\beta\alpha}(x,y)[\tau_{\alpha t}^b g^{\bk t}g_{b\bar
\beta}\p_{\bar k}(r^2)+\tau_{\bar \beta\bar t}^{\bar h} g^{\bar t
k}g_{\alpha\bar h}\p_k (r^2)]|_{(x_0,y)}
\end{split}
\end{equation}
where $(r^2)_{;\a\bbe}$ means the Hessian of $r^2(x,y)$ with respect to
the Riemannian connection,  $a,b,h,k,t\in\{1,\cdots,m+n\}$.  First of
all, we want to show that
\begin{equation}\label{eqn-laplacian2}
2g^{\bar\beta\alpha}(x,y)(r^2)_{;\a\bbe}(x_0,y)\leq
u(x,y)(\Delta^L_Mr^2)(x_0,y).
\end{equation}

Note that, by our choices of frames,
\begin{equation}\label{eqn-laplacian3}
2g^{\bar\beta\alpha}(x,y)(r^2)_{;\a\bbe}(x_0,y)=
2(r^2)_{;\a\ba}(x_0,y).
\end{equation}
By the assumption (3), the sectional curvature for the Levi-Civita
connection is nonpositive, we know $(r^2)_{;a\bar b}(x_0,y)$ is
positive definite, please see \cite{gw}, then for any fixed
$\alpha\in\{1,\cdots,m\}$,
\begin{equation}\label{eqn-laplacian4}
\begin{split}
2(D^2r^2)|_{(x_0,y)}(e_\a,e_\ba)
&\leq 2\textrm{trace}((D^2r^2)|_{(x_0,y)})g(e_\a,e_\ba)|_{(x_0,y)}\\
&= g_{\a\ba}(x_0,y)(\Delta^L_Mr^2)(x_0,y).
\end{split}
\end{equation}
Combining \eqref{eqn-laplacian3} and \eqref{eqn-laplacian4}, we can get
\eqref{eqn-laplacian2}. By Lemma 3.2 in \cite{vt} and  Lemma
\ref{lcle1} in this paper, under the assumptions of the curvature and
the torsion, we can get
\begin{equation}\label{eq-Laplacian r}(\Delta^L_Mr^2)(x_0,y)\leq
C(m,n,A)(1+r^2(x_0,y))^{\frac{\gamma+1}{2}}.
\end{equation}
Here $A$ is the same one as in the assumptions.
Submitting this to \eqref{eqn-laplacian2}, we can get
\begin{equation}\label{eqn-laplacian5}
\begin{split}
2g^{\bar\beta\alpha}(x,y)(r^2)_{;\a\bbe}(x_0,y)\leq
C(m,n,A)u(x,y)(1+r^2(x_0,y))^{\frac{\gamma+1}{2}}.
\end{split}
\end{equation}

Now we want to estimate the second term in the last equality of
\eqref{eqn-laplacian1}. Denoting $Q_{\alpha\bbe}=[\tau_{\alpha t}^b
g^{\bk t}g_{b\bar \beta}\p_{\bar k}(r^2)+\tau_{\bar \beta\bar t}^{\bar
h} g^{\bar t k}g_{\alpha\bar h}\p_k (r^2)]|_{(x_0,y)}$, it is a
2-tensor on $T^{(1,0)}_yY$. Choose a unitary basis
$\{s_{m+1},\cdots,s_{m+n}\}$ on $T^{(1,0)}_yY$, and extents it to
$\{s_1,\cdots,s_m,s_{m+1},\cdots,s_{m+n}\}$ as a unitary basis on
$T^{(1,0)}_{(x_0,y)}X\times Y$. Taking a vector $W=W^\alpha s_\alpha
\in T^{(1,0)}_yY$, we have
\noindent
\begin{equation*}
\begin{split}
|W^\alpha Q_{\alpha\bbe}
W^{\bbe}|&=|\vv<\tau(W,s_t),\overline{s_b}>\vv<s_b,\overline{W}>\overline{s_t}(r^2)||_{(x_0,y)}\\
&\leq |\vv<\tau(W,s_t),\overline{W}>\overline{s_t}(r^2)|_{(x_0,y)}\\
&\leq 2nA(1+r^2(x_0,y))^{\frac{\gamma+1}{2}}\|W\|^2_{g_{(x_0,y)}}.
\end{split}
\end{equation*}
Here we have used the  assumption (4) about the restriction on the
torsion with respect to the canonical connection.
So we get $$Q_{\alpha\bbe}\leq
2nA(1+r^2(x_0,y))^{\frac{\gamma+1}{2}}g_{\alpha\bbe}(x_0,y).$$
Hence

\noindent
\begin{equation*}\label{eqn-laplacian6}
\begin{split}
g^{\bbe \alpha}(x,y)Q_{\alpha\bbe}
&\leq 2nA(1+r^2(x_0,y))^{\frac{\gamma+1}{2}}g^{\bbe
\alpha}(x,y)g_{\alpha\bbe}(x_0,y)\\
&=2nAu(x,y)(1+r^2(x_0,y))^{\frac{\gamma+1}{2}}.
\end{split}
\end{equation*}
Combing this with \eqref{eqn-laplacian1} \eqref{eqn-laplacian5}, we
have
\begin{equation}\label{eqn-laplacian7}
\Delta f^2|_{(x,y)}\leq C_6 u(x,y)(1+r^2(x_0,y))^{\frac{\gamma+1}{2}}
\end{equation}
for some constant $C_6$ depends on $A, m, n$.

Let
\begin{equation}
w(x,y)=u(x,y)(C_0+r^2(x_0,y))^{-\delta}
\end{equation}
where $C_0\geq 1$ is a constant to be determined.
Noting that $0<\delta<\frac{1}{2}$, setting $v(x,y)=r^2(x_0,y)$, by
(3.5) in \cite{ty} or Corollary \ref{cor-Laplacian}, we have
\begin{equation*}
\begin{split}
\Delta w=&(C_0+v)^{-\delta}\Delta u-2\delta(C_0+v)^{-1-\delta}\vv<\nabla
u,\nabla v>\\
& -u\delta(C_0+v)^{-1-\delta}\Delta
v+u\delta(\delta+1)(C_0+v)^{-2-\delta}\|\nabla v\|^2\\
\geq &(C_0+v)^{-\delta}\Delta u-u\delta(C_0+v)^{-1-\delta}\Delta
v-2\delta(C_0+v)^{-1}\vv<\nabla w,\nabla v>\\
\geq &(C_0+v)^{-\delta}\Delta u-u\delta(C_0+v)^{-1-\delta}\Delta
v-2\delta(C_0+v)^{-1}\|\nabla w\|\cdot \|\nabla v\|.
\end{split}
\end{equation*}
Since $C_0\geq 1$, submitting the estimations of
\eqref{pt2lulb1}\eqref{eqn-gradient} and \eqref{eqn-laplacian7} to this
inequality, we can get
\begin{equation*}
\begin{split}
&\Delta w(x,y)\\
\geq&-2A(1+r^2(x,y))^{\gamma} w(x,y)+2B(1+v)^{-\delta}(C_0+v)^{\delta}
w^2(x,y)\\
&-4\delta|\nabla
w|(x,y)w^\frac{1}{2}(x,y)-C_6\delta(C_0+v)^{-1+\delta}(1+v)^\frac{\gamma+1}{2}w^2(x,y)\\
\geq&-2A(1+r^2(x,y))^{\gamma}w(x,y)+\Big(2B-C_6\delta(C_0+v)^{-\frac{1-2\delta-\gamma}{2}}\Big)w^2(x,y)\\
&-4\delta|\nabla w|(x,y)w^\frac{1}{2}(x,y).\\
\end{split}
\end{equation*}
Since $\gamma+2\delta<1$, we can choose $C_0$ large enough depending on
$C_6,\delta,\gamma,B$ such that
\begin{equation}
\Delta w\geq Bw^2-2A(1+r^2)^\gamma w-4\delta\|\nabla w\|w^\frac{1}{2}.
\end{equation}
So, a similar cut-off argument in the proof of Lemma \ref{lemls1} will
imply that
\begin{equation}\label{pt2ll2}
w\leq C_7(1+r^2)^\gamma
\end{equation}
where $C_7$ is positive constant depending on $A, B, m, n, \delta.$
\end{proof}

\begin{proof}[Proof of Theorem \ref{pm2}]

First of all, we may assume $M$ is simply connected because the
distance function in the universal cover of $M$ is no less than the
distance function of $M$. Suppose there is a complete metric $g$ on $M$
satisfying the assumptions of the theorem.
Let us choose the section $V$ as in the proof of Theorem \ref{pm1}, and
set $f(x)=|V|^2_{g(x,q)}$. By assumption (2) and the result of
\cite{k}, the holomorphic bisectional curvature of $X_q$ at $x$ is also
less than or equal to $-B(1+r^2(x,q))^{-\delta}$. By \eqref{t1-pf-v1}, we
have
\begin{equation}\label{pt2lp1}
\Delta_{X_q}f(x)\geq
2mB(1+r^2(x,q))^{-\delta} f(x).
\end{equation}
By Lemma \ref{t2-p-l1}, we can get
\begin{equation}\label{pt2lp2}
0<f(x)\leq C(1+r^2(x,q))^\delta.
\end{equation}

Now we want to show that
\begin{equation}\label{pt2lp3}
\Delta_{X_q}r^2(x,q)\leq C_8(1+r^2(x,q))^{\frac{1+\gamma}{2}}
\end{equation}
for some positive constant $C_8$ independent of $x$.
For any fixed point $x\in X_q$, choose an   holomorphic coordinate
$(z^1,\cdots,z^m)$ of $X_q$ at $x$ such that the induced metric $g^q$
on $X_q$ satisfies $g^q_{\alpha\bar\beta}(x)=\delta_{\alpha\beta},$
here $\alpha, \beta\in\{1,\cdots, m\}.$ Setting
$\varphi(x,y)=r^2(x,q)$, from Lemma \ref{lem-submfd}  and Corollary
\ref{cor-Laplacian}, we can get
\begin{equation}\label{pt2lp4}
\begin{split}
\Delta_{X_q}\varphi&=2\sum_{\alpha=1}^m(\varphi)_{\alpha\bar \alpha}\\
&=2\sum_{\alpha=1}^m(\varphi)_{;\alpha\bar \alpha}+[{\tau}_{\alpha
\beta}^\alpha  \p_{\bar \beta}(\varphi)+{\tau}_{\bar \alpha\bar
\beta}^{\bar \alpha} \p_\beta (\varphi)].
\end{split}
\end{equation}
Clearly
\begin{equation}\label{pt2lp4-2}
|{\tau}_{\alpha \beta}^\alpha  \p_{\bar \beta}(\varphi)+{\tau}_{\bar
\alpha\bar \beta}^{\bar \alpha} \p_\beta (\varphi)|\leq
C(1+\varphi)^{\frac{1+\gamma}{2}}
\end{equation}
for some constant $C$ independent of $x$.   Noting that the sectional
curvature for the Levi-Civita connection is nonpositive, by \cite{gw},
we can get $$2\sum_{\alpha=1}^m(\varphi)_{;\alpha\bar \alpha}\leq
\Delta^L_M \varphi.$$
By \eqref{eq-Laplacian r}, at $(x,q)$ we have $\Delta^L_M \varphi \leq
C(1+\varphi)^{\frac{1+\gamma}{2}}$ for some constant $C$ independent of
$x$. Combining this with \eqref{pt2lp4} and \eqref{pt2lp4-2}, we can
get \eqref{pt2lp3}.

Let $h(x)=\log f(x)-2\delta \log(C_9+r^2(x,q))$ where $C_9>1$ is some
constant.  Follow the proof of (3.10) in \cite{ty}, from \eqref{pt2lp2}
\eqref{pt2lp3} and the assumption $\gamma+2\delta<1$, we can get if $C_9$ is
big enough, then at a maximum point $(\bar x, q)\in X_q$
\begin{equation}\label{pt2lplll}
0\geq \Delta_{X_q}h(\bar x)>0.
\end{equation}
Hence we have a contradiction. Therefore Theorem \ref{pm2} holds.
\end{proof}

\end{document}